\documentclass[12pt,a4]{article}
\usepackage{mathptmx}       
\usepackage{helvet}         
\usepackage{courier}        
\usepackage{makeidx}         
\usepackage{graphicx}        
\usepackage{multicol}        
\usepackage[bottom]{footmisc}
\usepackage{amsmath,amssymb,amsfonts}        

\newtheorem{remark}{Remark}
\newtheorem{theorem}{Theorem}

\newtheorem{lemma}{Lemma}

\makeindex             

\newcommand{\exn}{{\bf E}}
\newcommand{\pr}{{\bf P}}
\newcommand{\N}{\mathbb{N}}
\newcommand{\R}{\mathbb{R}}

\newcommand{\al}{\alpha}
\newcommand{\be}{\beta}

\newcommand{\la}{\lambda}

\newcommand{\ga}{\gamma}
\newcommand{\si}{\sigma}
\newcommand{\de}{\delta}

\newcommand{\deq}{\stackrel{d}{=}}
\newcommand{\ep}{\varepsilon}

\newcommand{\ind}{{\bf 1}}

\newcommand{\ccdot}{{\, \centerdot}}

\title{On ruin probabilities in the presence of risky investments and random switching}
\author{Ying He$^1$  and Konstantin Borovkov$^2$ }
\date{}

\begin{document}

\footnotetext[1]{School of Mathematics and Statistics, The University of Melbourne, Parkville 3010, Australia; e-mail: hey4@student.unimelb.edu.au.}
\footnotetext[2]{School of Mathematics and Statistics, The University of Melbourne, Parkville 3010, Australia; e-mail: borovkov@unimelb.edu.au.}

%
%
\maketitle

\begin{abstract}\noindent 
	We study the asymptotic behavior of ruin probabilities, as the initial reserve goes to infinity,  for  a reserve process model where claims arrive according to a renewal process, while 	between the claim times the process has the dynamics of geometric Brownian
	motion-type It\^o processes with time-dependent random coefficients. These coefficients are ``reset" after each claim time, switching to new values independent of the past history of the process. We use the implicit renewal theory to obtain power-function bounds for
	the eventual ruin probability. In the special case  when the random drift and diffusion coefficients of the investment  returns process remain unchanged between consecutive claim  arrivals,  we obtain  conditions for existence of Lundberg's exponent for our model ensuring the power function behaviour for the ruin	 probability.



\end{abstract}

\section{Introduction and the main result}
\label{sec:1}

In the classical Cram\'er--Lundberg collective risk model (going back to a 1903 F.~Lundberg's work), the insurance company reserve process~$X$ is assumed to have dynamics of the form
\begin{align}
\label{CL}
X(t) = u + ct - \sum_{1\le j\le N(t)} \xi_j, \qquad t\ge 0,
\end{align}
where $c$ is a constant premium payment rate, $N$ is a Poisson process of claim epochs~$T_j,$ $j\ge 1,$  and $ \xi_j$, $j\ge 1,$ are  positive i.i.d.\ random variables modeling claim sizes made at the respective claim times, their sequence being independent of~$N.$

The main question posed in the context of this model was  on the behavior of the ultimate ruin probability
\[
\psi (u):= \pr \Bigl(\inf_{t>0} X (t) <0\Bigr)
\]
as  the initial reserve~$u$ tends to infinity. Clearly, in model~\eqref{CL}, the ruin  ($X$~turning negative) can only occur   at a claim time. Hence one deals here with a question on the asymptotic behavior of the distribution tail of the global maximum of a random walk with jumps of the form $Z_j:=\xi_j - c\tau_j,$ where $\tau_j:= T_j - T_{j-1},$ $j\ge 1$ (setting $T_0:=0$). Hence $\psi (u)<1$ for $u>0$ and $\psi (u) \to 0$ as $u\to\infty$ once the safety loading condition
\begin{align}
\exn Z \equiv \exn \xi - c\exn \tau <0
\label{SafLoaCon}
\end{align}
is met (here and in what follows, we use the convention that $Z \deq Z_1,$ $\xi \deq \xi_1,$ {\em etc.}). But what   can one say  about the rate at which $\psi (u)$ vanishes at infinity? The most famous result in this classical setting is the celebrated Cram\'er--Lundberg approximation that holds in the case of exponentially light tails and can be stated as follows.

For a random variable $V$,   denote by
\begin{align}
\phi_V(q):= \exn e^{qV} , \quad  q\in \R
\qquad
q_V:=\sup\{q\in \R:\phi_V(q)<\infty\}
\label{MGF_not}
\end{align}
its moment generating function and the right end-point of the interval on which the latter is finite, respectively.
If  $\phi_Z(q_Z)\ge 1$ then, under condition~\eqref{SafLoaCon}, there exists a unique solution $\ga >0$ to the equation $\phi_Z(q)=1$ and if, in addition, $\phi_Z(q_Z) > 1$ or $\phi_Z(q_Z) = 1$ and $\phi'_Z(q_Z- )<\infty$ then
\begin{align}
\psi (u) = C e^{-\ga u} (1+ o(1)) \quad \mbox{as\quad  $u\to\infty$}
\label{CL_approx}
\end{align}
where the constant~$C$ admits a closed-form expression
(see e.g.\ Section~22 in~\cite{Bo72} or  Section~I.4d in~\cite{AsAl10}). It turns out that approximation~\eqref{CL_approx} is rather sharp: there is an $\ep>0$ such that the remainder term~$o(1)$ in it can be replaced with $o(e^{-\ep u})$. Moreover, under the same moment assumptions on the distribution of~$Z$, approximation~\eqref{CL_approx} also holds for the  Sparre Andresen model that differs from~\eqref{CL} only in that the process~$N$ is just a renewal process (so that the   inter-claim times~$\tau_j$ are general  positive i.i.d.\ random variables). In this case, the remainder term will be decaying exponentially fast  under the additional assumption that the distribution of~$Z$ contains an absolutely continuous component (p.~129 in~\cite{Bo72}).

Note that in the case where $\phi_Z(q_Z) = 1$ and $\phi'_Z(q_Z- )=\infty$, the problem on the  asymptotics of~$\psi (u)$ is more difficult and the asymptotic behavior of this probability as $u\to \infty$ can have a different form, see e.g.\ p.\,136 in~\cite{Bo72}  and Section~6.5 in~\cite{BoBo08}.

Of course, the Cram\'er--Lundeberg model~\eqref{CL} and its Sparre Andersen extension are oversimplifications of real-life situations. These models assume that all the reserves of the insurance company are kept in a safe bank account. Over the last two decades, several authors turned their  attention to more realistic models in which the reserve capital can be invested in a risky financial asset (considering a single  risky asset is reasonable due to the common practice of investing in a market portfolio or an index). Models with surplus generating process and investments in risky asset modelled by L\'evy processes were discussed, e.g., in~\cite{Pa93, Pa98, PaGj97}. In particular, it was noted in~\cite{PaGj97} that the ultimate ruin probability and the Laplace transform of the ruin time are solutions to suitable  boundary value problems for the respective   integro-differential equations.

A discrete time model with stochastic interest rates and returns was considered in~\cite{Ny99}, the main results (obtained using the ``crude" large deviation theory)   included power asymptotic behavior of the ruin probability as a function of the initial reserve.
A power function ruin probability asymptotics behavior was also obtained in~\cite{Pa02} for the L\'evy processes-based models under suitable conditions, basing on the results from~\cite{Pa98}.  Assuming~\eqref{CL} (and also allowing a more general L\'evy process model) and that the risky investment returns follow an independent geometric L\'evy process, power function bounds for~$\psi (u) $ were obtained in~\cite{KaNo02}. A power function asymptotic behavior was   obtained in~\cite{FrKaPe02} for a modification of the classical  model~\eqref{CL} with investments in a risky asset with price following an independent geometric Brownian motion (BM) process with mean return~$\mu\in \R$ and volatility~$\si>0$ (as in~\eqref{PeZeMod} below, but with a constant~$c(s)\equiv c$). 
Assuming that the claim sizes are exponentially distributed and setting $\be:=2\mu /\si^2 -1 $, it was shown in~\cite{FrKaPe02} that
\begin{align}
\psi (u)   = C u^{ -\be}(1+o(1)) \quad \mbox{as\quad  $u\to \infty$}
\label{power_ass}
\end{align}
for some constant $C>0$ when $\be>0$ (and that  $\psi (u)\equiv 1$ when $\be<0$). For claims with a general distribution such that $\exn \xi_1^\beta <\infty$, were obtained upper and lover power bounds with the right-hand sides of the form $Cu^{ -\be}$ for come constants~$C$.  

Note that the presence of the moment condition on~$\xi_1$ (here and in our Theorem~\ref{T1} below) is quite natural as for heavy-tailed claim distributions, the asymptotics of the ruin will be governed by the distribution tail of the ``integrated tail law"  for~$\xi_1$ when that tail dominates~$u^{-\be}$ (cf.~\cite{AlCoTh12} and Chapter~X in~\cite{AsAl10}).

These results were  extended in~\cite{PeZe06} to a modification of the above model with a variable premium payment rate $c(t)$  yielding the following dynamics:
\begin{equation}
\label{PeZeMod}
X(t) = u + \int_0^t c(s)\,ds
+ \int_0^t \mu X(s)\, ds
+ \int_0^t \si X(s)\, dW(s)
- \sum_{j\le N(t)}\xi_j,
\end{equation}
where  $W$ is a BM process independent of $N$ and $\{\xi_j\}$, the coefficients~$\mu$ and~$\si$ are constant, and  $c (t)=c(t,X)\in [0,\overline{c}]$ (with a constant $\overline{c}<\infty$) is a bounded adapted function such that there exists a unique strong solution to the above equation. Upper and lover bounds with the right-hand sides of the form $Cu^{ -\be} $ were obtained under   appropriate moment conditions on~$\xi_1$, whereas exact asymptotics of the form~\eqref{power_ass} were established for  generally distributed~$\xi_1$ (satisfying $\exn \xi_1^{\be +\de}<\infty $ for some $\de>0$) in the special case where $c(t)=c_1 e^{\ga t}$ for some $\ga\le 0.$ The toolbox used in that paper, as in some other previous work as well, was based on the implicit renewal theory.

It may seem paradoxical at  the first glance that, in all these papers establishing power asymptotics of the form~\eqref{power_ass}, the distributions of the ``main source of risk"\,---\,the claims made against the insurer\,---\,could have  a finite exponential moment, as in the case leading to the much faster exponential decay~\eqref{CL_approx}. This means that investing in a risky asset (even with significant mean positive returns) dramatically increases the riskiness of the insurance business. In Remark~\ref{R5} below we will provide an intuitive explanation of the emergence of the power behavior at infinity for~$\psi$. Roughly speaking, it is due to the closeness of the dynamics of an embedded discrete time process (the values of the risk process~$u-X$ at the claims times) to those of the exponential of a random walk with i.i.d.\ jumps and negative trend. The ruin occurs when the global supremum of that walk is ``large", of the order of magnitude of~$\ln u$, and the probability of this  has the form of the right-hand side of~\eqref{CL_approx}, with $u$ replaced by~$\ln u$

Over the last few years, several authors turned their attention to versions of model~\eqref{PeZeMod} with random switching. In~\cite{ElKa20}, it was assumed that the geometric BM process modelling the dynamics of the risky asset has stochastic drift and volatility coefficients: $\mu = \mu_{\theta (t)},$  $\si = \sigma_{\theta (t)},$ where $\{\theta (t)\}_{t\ge 0}$ is a time-homogeneous (hidden) Markov chain with state space $\{0,1\}$ independent of all the other stochastic ingredients  of the model. Using implicit renewal theory, the authors derived two-sided power function bounds of the form
\begin{align}
0<\liminf_{u\to\infty } u^\beta \psi (u)
\le \limsup_{u\to\infty } u^\beta \psi (u) <\infty
\label{2sided}
 \end{align}
for the ruin probability. These results were extended in~\cite{KaPe22} to the case where $\{\theta (t)\}_{t\ge 0}$  has an arbitrary finite state space.

In~\cite{EbKaTh22} a Sparre Andersen type model was considered,  where the dynamics of the risky asset used for investment was given by a general L\'evy process $\{R(t)\}_{t\ge 0}$ (with the assumption that its jumps are always greater than $-1$):
\[
X (t)= u+\int_0^t X (s-)dR (s)- \sum_{j\le N(t)}\xi_j, 
\]
where $\{N(t)\}_{t\ge 0}$ is now a renewal process (all the components of the model were, as usual, assumed to be independent). Using recent results from  the  theory of distributional equations, the authors   derived for this model  two-sided power function bounds of the form~\eqref{2sided}.

In the present note, we   extend~\eqref{PeZeMod} to another version of the Sparre Andersen-type model with investment in a risky asset that  involves random switching. To formally describe our model, in addition  to the i.i.d.\ sequence  $\{\xi_\ccdot\}$ of claim sizes (as above), introduce an independent of it i.i.d.\ sequence of quadruples 	
\begin{align}
\label{quadr}
(\mu_n (\cdot), \si_n  (\cdot), \tau_n, W_n (\cdot)), \quad n\ge 1,
\end{align}
and (independent) filtrations $\{ \mathbb{H}_n=\{\mathcal{H}_n (t), t\ge 0\}\}_{n\ge 1}$, where   $W_n$ is a standard Wiener process  which is a martingale w.r.t.\ filtration~$\mathbb{H}_n$, while the process   $\mu_n  $  is adapted to $\mathbb{H}_n$ and locally integrable a.s., $ \si_n $ is progressively measurable (w.r.t.~$\mathbb{H}_n$) and locally square-integrable a.s., and $\tau_n>0$ are stopping times w.r.t.~$\mathbb{H}_n$ (in particular, they may be independent of~$W_n$, assuming~$\mathbb{H}_n$ large enough). About $c(t)$ we will   assume, as in~\cite{PeZe06}, that it is right-continuous and takes values in $[0,\overline{c}]$ with some $0< \overline{c}<\infty$ and is adapted in an appropriate way (omitting technical details to avoid making exposition too cumbersome) such that there exist unique strong solutions to the equations describing our model.

Our reserve process follows the dynamics of~\eqref{PeZeMod}, where the drift and diffusion coefficients~$\mu$ and~$\si$ are  random   processes of the form
\begin{align*}
\mu (t)=\sum_{n=1}^\infty \mu_n (t-T_{n-1})\ind_{[T_{n-1}, T_n)}(t), \quad
\si (t)=\sum_{n=1}^\infty \si_n (t-T_{n-1}) \ind_{[T_{n-1}, T_n)}(t),
\end{align*}
while  $N(t)=\sum_{n=1}^\infty \ind_{(0,t]}(T_n),$ $T_n:=\sum_{i=1}^n\tau_i$,  is the renewal process generated by  the inter-arrival times $\tau_n>0$. We assume that
\[
\mu_n(t)\ge   \underline{\mu}>-\infty, \quad
0<\si_n(t)\le\overline{\si}<\infty\quad{\rm a.s.}
\]	
for some constant $\underline{\mu}, \overline{\si}$.

Thus, according to the suggested  model, our insurance company  commences at time $t=0$ with an initial endowment~$u$, faces a  renewal-reward claims process with  claim sizes $\xi_n$ and inter-claim times $\tau_n$, and  receives premium inflow at a bounded non-negative random rate~$ c( t).$ During the time period $(T_{n-1}, T_n)$,  the   company obtains a rate of return following  a diffusion process  with random time-dependent drift coefficient~$\mu_n  $ and volatility~$\sigma_n$, which are ``switched'' to~$\mu_{n+1}$ and~$\si_{n+1}$ at time~$T_n.$ The random regime switching for the investment  component may be related to changing the investment policy or insurer's economic environment  following claim payments. Considering the proposed model is also suggested  by the inner logic of the mathematical problem per se.

To state our main results, we first need to introduce some notations. Following the standard approach used, in particular, in~\cite{KaNo02} and~\cite{PeZe06}, we note that ruin for this model can only occur at one of the claim times~$T_n$. Therefore,  for the ruin probability analysis,  it suffices to consider the  embedded discrete time process $\{S_n:=X (T_n)\}_{n\ge 0}$ (setting $T_0:=0$) since
\begin{align}
\label{psi_S}
\psi (u) = \pr \bigl(\inf_{n\ge 1} S_n<0\bigr).
\end{align}

The dynamics of~\eqref{PeZeMod} inside intervals $[T_{n-1}, T_n)$ are those of solutions to linear stochastic differential equations with the respective initial values~$S_{n-1}$. Using the available in closed form solutions to such problems  (see e.g.\ Chapter~9 in~\cite{Ma13}), noting that $S_n = X(T_n-)-\xi_n,$ and introducing notations
\[
K_n (s):= \int_s^{\tau_n} (\mu_n (u)- \si_n^2(u)/2)du, \quad
Z_n (s):= \int_s^{\tau_n}  \si_n (u) d W_n(u), \quad   s\in [0,\tau_n],
\]
$K_n:= K_n (0),  $ $Z_n:= Z_n(0) $,
\begin{align}
\label{zeta}
\nu_n: =-K_n-Z_n,\quad  \lambda_n:= e^{- \nu_n}, \quad
\zeta_n:=\int_0^{\tau_n} e^{K_n (s)+Z_n(s)} c(T_{n-1}+s)ds - \xi_n,
\end{align}
we obtain that
\begin{align}
\label{S_n}
S_n= \la_n S_{n-1} +\zeta_n, \quad n\ge 1, \quad S_0 = u,
\end{align}
Note that, due to our assumptions,  $\{(K_n, \nu_n)\}_{n\ge 1}$   is an i.i.d.\ sequence, whereas
$\{\zeta_n\}_{n\ge 1}$   does not need to be so.

Recall that,  for sequences of random elements,  we agreed to omit for brevity's sake  the subscript~$n$ in the case where $n=1$.

Referring to~\eqref{MGF_not}, we will use the following lemma to introduce one more notation.

\begin{lemma}\label{L1}
	If $\phi_\nu (q_\nu)\in (1,\infty],$ $\exn \tau <\infty$ and $\exn K\in (0,\infty)$ then $q_\nu>0$ and there exists a $\beta\in (0,q_\nu)$ such that $\phi_\nu (\beta) =1.$
\end{lemma}

We will refer to $\beta$ from Lemma~\ref{L1} as the {\em Lundberg exponent\/} for our model.

\begin{remark}\label{R1}\rm 
	 Note that as  $\phi_\nu$ is left-continuous on $(0,q_\nu)$, one has $\be <q_\nu. $ Therefore  $\exn \la^{-(\be +\de)}\equiv \phi_\nu (\be +\de)<\infty$ for any $\de\in (0,q_\nu-\be)\neq \varnothing$.
\end{remark}

Our   main result is stated in the following theorem.

\begin{theorem}\label{T1}
Assume that  $\phi_\nu (q_\nu)\in (1,\infty], $ $\exn \tau <\infty,$  $\exn K\in (0,\infty),$ and $\exn \xi^{\be  }<\infty$ for some  $\de >0$, where    $\beta $ is the Lundberg exponent for our model.
Then
	\begin{align}
	\label{5.7}
	\limsup_{u\to\infty} u^\be \psi (u) & \le C_+.
	\end{align}
If, in addition,  $(\mu(\cdot), \si (\cdot), \tau)$ and $W(\cdot)$ are independent,  $\exn \xi^{\be	+\de}<\infty$ for some  $\de >0$ and
\begin{align}
\label{cond_tau}
q_\tau >
 \be  ^2\overline{\si}^2 /2+ \be (\overline{\si}^2/2 -\underline{\mu})^+	
	\end{align}	
then
\begin{align}
	\label{5.8}
	\liminf_{u\to\infty} u^\be \psi (u) & \ge C_- .
\end{align}
Here $0<C_-\le C_+ <\infty $ are some constants.
\end{theorem}

\begin{remark}\label{R2}\rm 
The existence of the Lundberg exponent $\be>0$ is ensured  since  the  conditions of  Lemma~\ref{L1} are clearly met under the assumptions of Theorem~\ref{T1}. Without loss of generality,  in what follows we will assume about   the~$\de$ from  the conditions of Theorem~\ref{T1} that $\de\in (0, q_\nu-\be)$ (see Remark~\ref{R1}).
\end{remark}

\begin{remark}\label{R2a}\rm 
Condition $\exn K= \exn \int_0^\tau (\mu (u) -\si^2 (u)/2) du>0$ means that ``volatility'' $\si (t)$ cannot be ``large"   in some average sense. Recall that $\si^2/2>\mu$ implies certain ruin in the models with constant $\mu$ and $\si$ considered in~\cite{FrKaPe02} and~\cite{PeZe06}.
\end{remark}

\begin{remark}\label{R3}\rm 
Observe  that if $\phi_\tau (q)<\infty$ for any $q>0$ then condition~\eqref{cond_tau} is clearly  superfluous.
\end{remark}

The proof of Theroem~\ref{T1} is given in Section~\ref{S2}.

The existence of the Lundberg exponent~$\be$ is the key factor for establishing the power behaviour of the ruin probability.  Given the structure of our random variable~$\nu$, verifying the existence of such a~$\be$  in the general case is a complicated task. In Section~\ref{S3}, we will establish  a sufficient condition for the existence of the Lundberg exponent in the more tractable special case when
\begin{align}
\label{musi}
\mu_n (t)\equiv \mu_n(0)=:\mu_n \quad\mbox{and}\quad \si_n (t)\equiv \si_n (0)=:\si_n
\end{align}
do not depend on time (so that the random drift and diffusions coefficients for the return on investments process remain unchanged during each of the intervals $[T_{n-1},T_n),$ $n\ge 1$). Moreover, we will assume  that the components  $(\mu_n,\si_n)$, $\tau_n$ and $W_n(\cdot)$ of  our quadruples~\eqref{quadr} are jointly independent. The problem admits in this case an elegant  solution: it turns out that the answer  (given in Theorem~\ref{T2} stated and proved in Section~\ref{S3}) basically depends on ``concentration of probability" in vicinity of a certain straight line tangent to the support of the distribution of the random vector~$(\mu, \si^2/2)$.

\section{Proof of Theorem~\ref{T1}}\label{S2}

\noindent{\em Proof of Lemma~\ref{L1}.}  That $q_\nu > 0$ in clear since  $\phi_{\nu} (q_{\nu})>1.$ Further,   
\begin{align}
\label{e_nu}
\exn\nu =- \exn K -\exn Z=- \exn K<0
\end{align}
as $\exn Z=\exn  \int_0^{\tau }  \si  (u) d W (u)=0 $ by the optional stopping theorem (note that $\exn|\nu|<\infty$).  Since~$\phi_\nu$ is a convex function and $\phi_\nu (q_\nu)\in (1,\infty],$ the existence of the claimed $\beta$ is equivalent to having  $ \phi'_\nu (0+)<0 ,$ which is an immediate consequence of~\eqref{e_nu}.\hfill   $\Box$

\medskip

\noindent{\em Proof of Theorem~\ref{T1}.} Our line of argument follows the overall logic employed in~\cite{PeZe06}. Iterating~\eqref{S_n} and setting $\Lambda_n:=\prod_{k=1}^n \la_k,$ $k\ge 1,$ we get
\begin{align}
\label{S_n_1}
S_n= \Lambda_nu + \Lambda_n \sum_{k=1}^n \Lambda_k^{-1} \zeta_k,\quad n\ge 1.
\end{align}

First we will prove the upper bound~\eqref{5.7}. Clearly,
\[
(Q_k, M_k):= \{(\xi_k,1)/\la_k\}_{k\ge 1}
\]
is an i.i.d.\ sequence. Set
\begin{align}
\label{R_n}
R_n :=  \sum_{k=1}^n Q_k \prod_{i=1}^{k-1}M_i, \quad n\ge 1,
\end{align}
with the usual convention that $\prod_{i=j}^k=1$ when $j>k.$ Since $Q,M>0,$ the sequence   $\{ R_n \}_{n\ge 1}$ is clearly  increasing so that
\begin{align}
\label{R}
 R_n\uparrow R\quad\rm a.s.
\end{align}
for some (possibly improper) random variable $R\le \infty.$

In view of~\eqref{zeta}, one has  $\zeta_k\ge -\xi_k$, $k\ge 1,$ and hence we obtain from~\eqref{S_n_1} that
\begin{align}
S_n \ge  \Lambda_n u - \Lambda_n \sum_{k=1}^n \Lambda_k^{-1} \xi_k
 =  \Lambda_n (u-R_n)\ge \Lambda_n (u-R ),\quad n\ge 1.
\end{align}
Hence it follows from~\eqref{psi_S} that
\begin{align}
\label{upper_R}
 \psi (u) \le \pr (R>u),\quad u>0.
\end{align}

\begin{remark}\label{R5}\rm 
One can clarify the emergence of the power decay for~$\psi$ as follows. Clearly, $U_n:=\sum_{k=1}^n \nu_k, $ $n\ge 1,$ is a random walk with i.i.d.\ jumps $\nu_k$ with  negative trend (see~\eqref{e_nu}) and $\phi_\nu (\be )=1$. Hence by the classical Cram\'er--Lundberg result~\eqref{CL_approx} for $\overline{U} := \max_{n\ge 1} U_n$, one has  $\pr (\overline{U} >w )\sim C e^{-\beta w}$ as $w\to\infty$.

Now in view of~\eqref{S_n_1}, ruin is equivalent to the event  $\Big\{\sup_{n\ge 1 }\sum_{k=1}^n (-\zeta_k)e^{U_k}>u\Big\}$ which actually occurs ``due" to a few terms in these sums, with~$k$ close to the point~$n'$ such that $\overline{U} =U_{n'}$ (cf.\ the argument in the proof of Theorem~4 in~\cite{BoVa06}). So one can expect that the probability of ruin behaves like $\pr (\overline{U} >\ln u)\sim C e^{-\beta \ln u}=Cu^{-\be}$ as $u\to\infty.$
\end{remark}

That $R$ is a proper random variable follows immediately from the following lemma, which is a direct consequence of Theomre~1.6 in~\cite{Ve79}:

\begin{lemma}\label{L2}
Let $\{(A_n, B_n)\}_{n\ge 1}$ be an i.i.d.\ sequence of bivriate random vectors, and
\begin{align}
 \label{Z_n}
	Z_n (x):= x\prod_{j=1}^n A_j + \sum_{k=1}^n B_k \prod_{j=1}^{k-1} A_j, \quad n\ge 1,\quad x\in \R.
\end{align}
Assume that  $\exn \ln |A|<0$ and $\exn (\ln |B|)^+<0,$  where $z^+:= \max\{0,z\},$ $z\in \R. $ Then  $Z_n (x)\to Z$ in distribution as $n\to\infty$   for all $x\in \R,$ where the distribution of the proper random variable~$Z$ satisfies the random equation
\begin{align}
\label{eq_Z}
Z\deq B+AZ,
\end{align}
$(A,B)$ and $Z$ on the right-had side being  independent of each other.
\end{lemma}	

Indeed, our sequence~\eqref{R_n} is of the form~\eqref{Z_n} with $x=0$ and $(A_n, B_n)=(M_n, Q_n),$ $n\ge 1,$ and $ \exn \ln |A|=\exn \ln \la^{-1}=\exn \nu<0$ by~\eqref{e_nu},  whereas
\[
\exn (\ln |B|)^+
  = \exn (\ln  (\xi/\la))^+
  = \exn (\ln   \xi +\nu)^+
  \le \exn (\ln   \xi)^+  +\exn \nu ^+
<\infty
\]
as $\exn \xi^\be <\infty$ and $\exn|\nu|<\infty$ (cf.~Lemma~\ref{L1}).

Hence, by Lemma~\ref{L2}, the sequence  $ \{R_n\}$ converges as $n\to\infty$ in distribution to a proper random variable, which implies that the a.s.\ limit~$R$ from~\eqref{R} is proper as well and satisfies the random equation
\begin{align}
\label{eq_R}
R\deq Q+ MR,
\end{align}
where $R$ and $(M,Q)$ on the right-hand side are independent of each other.

Now to complete the derivation of the desired upper bound using~\eqref{upper_R} it remains to turn to the implicit renewal theory.  We will make use of the following lemma which is a direct consequence of Theorem~4.1 in~\cite{Go91}.

\begin{lemma}\label{L3}
Assume that the distribution of a bivariate random vector $(A,B)$ with $A\ge 0$ a.s.\ is such that, for some $\al>0$, 	
\[
\exn A^\al =1, \quad \exn A^\al (\ln A)^+<\infty, \quad \exn |B|^\al <\infty,
\]
while the conditional distribution of $\ln A$ given $A\neq 0$ is non-arithmetic. Then solution to~\eqref{eq_Z} satisfies
\[
\lim_{u\to\infty} u^\al \pr (Z>u)=C,
\]
where $C:= \exn \big[ ((B+AZ)^+)^\al- (AZ^+)^\al\big]    /(\al \exn A^\al \ln A) \in (0,\infty).$
\end{lemma}

To apply this lemma to our equation~\eqref{eq_R} with $(A,B)=(M,Q)$ and $\al =\be$, it suffices to note that $\exn M^\be = \phi_\nu (\be)=1$, $\exn Q^\be =\exn \xi^\be M^\be =\exn \xi^\be \exn M^\be =\exn \xi^\be<\infty$ due to independence, and  $\exn M^\be (\ln M)^+ = \exn e^{\be \nu} \nu^+  <\infty$ since $\exn e^{(\be+\de) \nu} <\infty$ for some $\de >0$ (see Remark~\ref{R2}). That $\ln M$ given $M\neq 0$ is non-arithmetic is obvious from the definition of $M=e^\nu$ and the presence of the It\^ o integral in~$\nu$. This completes the proof of the upper bound~\eqref{5.7}.

Now we will proceed  to proving the lower bound~\eqref{5.8}. The main tool here is the following assertion  from~\cite{Le86} (see also~\cite{Go91} and~\cite{Ny01}).

\begin{lemma}
\label{L4}
 	Assume that $Y$ satisfies the equation
\begin{align}
\label{5.3}
 	Y\deq B+AY^+,
\end{align}
where $(A,B)$ and $Y$ on the right-hand side are independent of each other, $A>0$ a.s., and the distribution of~$(A,B)$ is such that $\pr (A>1, B>0)>0.$ If, for some $\al,\de>0$, 		
\[
\exn A^\al =1, \quad \exn A^{\al+\de}  <\infty, \quad \exn |B|^{\al+\de}  <\infty,
\]
and $\ln A$ is absolutely continuous, then
\[
\lim_{u\to\infty} u^\al \pr (Y>u)=C+ o(u^{-h})
\]
for some positive constants $C$ and $h$. 	
\end{lemma} 	

To apply this result, we turn to representation~\eqref{S_n_1} and use the natural upper bound for~$\zeta_n:$
\begin{align*}
\label{barZ}
\zeta_n\le \overline{\zeta}_n:= \overline{c}\int_0^{\tau_n} \exp\{K_n (s) +Z_n(s)\}ds - \xi_n
\end{align*}
to get the inequality
\[
S_n\le \overline{S}_n:= \Lambda_nu + \Lambda_n \sum_{k=1}^n \Lambda_k^{-1} \overline{\zeta}_k
 = \Lambda_n (u-\overline{R}_n ),\quad n\ge 1,
\]
where
\begin{align*}
\overline{R}_n :=  \sum_{k=1}^n \overline{Q}_k \prod_{i=1}^{k-1}M_i, \quad
  \overline{Q}_n:=- \overline{\zeta}_n/\la_n, \quad n\ge 1.
\end{align*}
In view of~\eqref{psi_S}, this  implies the bound
\begin{align*}
\psi (u) \ge \pr \big(\inf_{n\ge 1 }\overline{S}_n <0\big)
 \ge \pr \big(\overline{R} >u), \quad \mbox{where} \ \ \overline{R}:=\sup_{n\ge 1} \overline{R}_n.
\end{align*}
Next we note that, since $\overline{R}_1=\overline{Q}_1$ and  $M_1>0,$ one has
\begin{align*}
\overline{R} &= \overline{Q}_1\vee \sup_{n\ge 2}\overline{R}_n
  =\overline{Q}_1\vee \Big(\overline{Q}_1 +M_1\sup_{n\ge 2}\sum_{k=2}^n \overline{Q}_k \prod_{i=2}^{k-1}M_i \Big)
\\
& = \overline{Q}_1\vee (\overline{Q}_1 +M_1\overline{R}' )
 = \overline{Q}_1 +M_1(\overline{R}'  )^+,
\end{align*}
where $\overline{R}':=\sup_{n\ge 2}\sum_{k=2}^n \overline{Q}_k \prod_{i=2}^{k-1}M_i \deq \overline{R}$ is independent of $(M_1, Q_1).$ Therefore our~$\overline{R}$ satisfies the random equation
\[
\overline{R} \deq \overline{Q}  +M ( \overline{R}  )^+,
\]
where $(M,Q)$ and $\overline{R}$ on the right-hand side are independent of each other. This relation  is exactly   of the form~\eqref{5.3}, and we will now verify whether the conditions of Lemma~\ref{L4} are met when   $(A,B)= (M,\overline{Q}),$  $\al=\be.$

First of all, it follows from Proposition~6.1 in~\cite{Go91} that   $\overline{R}$ is a proper random variable provided that $\exn \ln (1\vee \overline{Q})<\infty$. The latter  will immediately  follow from the condition $\exn |\overline{Q}|^{\be +\de}<\infty $ of Lemma~\ref{L4} that we need to verify. To demonstrate the latter relation, note that
\begin{align}
\overline{Q}& = -\frac{\overline{\zeta}}{\la}
 = \frac{ {\xi}}{\la} - \overline{c}e^{-K(0) -Z(0)} \int_0^\tau e^{K (s)+Z (s)} ds
\notag \\
& = \frac{ {\xi}}{\la} - \overline{c} \int_0^\tau \exp\bigg\{
- \int_0^s (\mu  (u)- \si ^2(u)/2)du  -\int_0^s  \si (u) d W (u)  \bigg\} ds.
\label{for_Q}
\end{align}
It is obvious  from  the elementary inequality $|x+y|^p\le (1\vee 2^{p-1})(|x|^p+|y|^p) ,$ $x,y,p>0,$ that it suffices to show that the absolute moments of the order $\be +\de$ are finite for both terms on the right-hand side. By independence, one has
\[
\exn \bigg|\frac{ {\xi}}{\la}\bigg|^{\be +\de}
 =\exn  {\xi}^{\be +\de}\exn\la ^{-(\be+\de)} =\exn  {\xi}^{\be +\de}\phi_\nu (\be+\de)<\infty
\]
in view of  Remark~\ref{R2}.

Next note that, due to our assumption about independence of $(\mu (\cdot),\si (\cdot), \tau)$  and $W(\cdot),$ one has
$
\big\{ -\int_0^s \si (u) dW(u)\}_{s\ge 0} \deq \{W(\Sigma (s))\big\}_{s\ge 0} ,
$
where we set   $\Sigma (s):= \int_0^s \si^2 (u)du, $ $s\ge 0.$ Therefore, putting $\overline{W}(t):=\max_{0\le s\le t} W(s), $ $t\ge 0,$ we get
\[
\max_{0\le s\le \tau} \biggl( -   \int_0^s  \si (u) d W (u)   \bigg)
\deq \max_{0\le s\le \tau} W(\Sigma (s))=\overline{W}(\Sigma (\tau))
\le \overline{W}(\overline{\si}^2  \tau) .
\]
Now, setting $\kappa (u):=\si ^2(u)/2-  \mu  (u)  $ and noting that  $\kappa (u)\le \overline{\kappa}:= \overline{\si} ^2/2 -\underline{\mu } $ a.s.,  we get for the second term on the right-hand side of~\eqref{for_Q} that
\begin{align}
\exn \bigg(\int_0^\tau \cdots ds\bigg)^{\be + \de}
 &\le \exn \bigg(e^{\overline{W}(\overline{\si}^2  \tau)}
 \int_0^\tau  e^{ \overline{\kappa} s} ds\bigg)^{\be + \de}
  \le   \exn \bigg(e^{\overline{W}(\overline{\si}^2  \tau) }\tau
  e^{ \overline{\kappa}^+\tau }\bigg)^{\be + \de}.
  \label{int_bd}
\end{align}
Due to the reflection principle, for any $a, t>0,$
\[
\exn e^{a \overline{W}(t)}=  2\exn (e^{a  {W} (t)}; W (t) >0) <  2\exn  e^{a {W} (t)}=2e^{a^2  t/2},\quad t>0,
\]
so conditioning the last expectation in~\eqref{int_bd} on $\tau$ and using independence,  we obtain that it is less than
\[
2 \exn e^{((\be + \de)^2\overline{\si}^2 /2+ (\be + \de)\overline{\kappa}^+)\tau  } \tau^{\be + \de}<\infty ,
\]
using assumption~\eqref{cond_tau} and choosing~$\de>0$ small enough. Thus we showed that $\exn|\overline{Q}|^{\be+\de}<\infty,$ which implies, in particular, that $\overline{R}$ is proper.

To verify the remaining assumptions of Lemma~\ref{L4}, we observe  that condition $\exn M^\be=1$ is met by Lemma~\ref{L1} and that $\exn M^{\be+\de}<\infty$ as explained in Remark~\ref{R2}. That  $M=  e^\nu $  is absolutely continuous follows from the presence of the It\^o integral in~$\nu$ and independence of~$W$ from the other participating random quantities. Thus it only remains to verify that $\pr (M>1,\overline{Q}>0)>0$. Setting
\[
V (t):=\int_0^{t} (\mu  (u)- \si ^2(u)/2)du+ \int_0^{t}  \si  (u) d W (u), \quad   t\ge 0,
\]
and choosing $a>0$ such that $b:= \pr ( \xi >a\overline{c})>0$, the previous  probability is clearly equal to
\begin{align*}
\pr (V(\tau) <0, \overline{\zeta}<0 )
& \ge  \pr \bigg(V(\tau)<0, \overline{c} \int_0^\tau e^{V(\tau) -V(s)}ds <\xi , a \overline{c} <\xi\bigg)
\\
& \ge b \pr \bigg(V(\tau) <0,  \int_0^\tau e^{V(\tau) -V(s)}ds < a  \bigg)
\\
& \ge b \pr \bigg(V(\tau) <0,  \int_0^\tau e^{  -V(s)}ds < a  \bigg)
\\
&  \ge b \pr \big(  \tau e^{  -\underline{V} (\tau)} < a \,|\, V(\tau) <0  \big)\pr(V(\tau) <0),
\end{align*}
where we put $\underline{V} (t):=\inf_{0\le s\le t}{V}(s)$. Obviously, $\pr(V(\tau) <0)>0$, and as $-\underline{V}(\tau) >0$ on the event $\{V(\tau) <0\}$ while~$a$ can be chosen arbitrary small, the product in the last line of the displayed formula is positive, establishing that the last condition of Lemma~\ref{L4} is met as well. This completes the proof of Theorem~\ref{T1}.\hfill{$\Box$}\medskip

\section{Lundberg's exponent when coefficients  $\mu_n(t)$ and $\si_n(t)$ do not depend on time}\label{S3}

In this section we assume satisfied condition~\eqref{musi} and also that the components $(\mu_n,\si_n)$, $\tau_n$ and $W_n(\cdot)$   of our quadruples~\eqref{quadr} are jointly independent. Under these assumptions, one has $\nu = -(\mu - \si^2/2) - \si W(\tau)$. Introducing the random vector $\Theta:= (\mu,\si^2/2),$ setting 
\[
u (q):= (-q, q (q+1)),\quad q\in \R, 
\] 
and conditioning,  we get
\begin{align}
\phi_\nu (q) & = \exn e^{q (- (\mu-\si^2/2)\tau - \si W(\tau))}
=\exn e^{ -q(\mu-\si^2/2)\tau +  q^2 \si^2 \tau/2}
\notag
\\
& =\exn \phi_\tau (-q(\mu-\si^2/2)  +  q^2 \si^2 /2)
= \exn \phi_\tau (\langle u(q), \Theta \rangle),
\label{phi_now}
\end{align}
where $\langle\cdot, \cdot \rangle$ stands for the inner product in~$\R^2.$

Note that our  key condition $\exn K>0$ for the existence of the Lundberg exponent is equivalent in the case under consideration to
\begin{align}
\label{mu_si_0}
\exn (\mu -\si^2/2)>0
\end{align}
(assuming that $\exn \tau<\infty$), which is a  ``mean version" of the condition $2\mu/\si^2 -1>0$ under which the asymptotics~\eqref{power_ass} was established in the case of constant deterministic $\mu$ and~$\si$ in~\cite{PeZe06}.

Assuming that  the above condition is met, the case $q_\tau =\infty$ is trivial: it is clear from  Lemma~\ref{L1} and~\eqref{phi_now} that $\beta$ will then always exist. So we will only consider the case where
\begin{align}
\label{tau_infi}
 q_\tau <\infty,\quad   \phi_\tau (q_\tau)=\infty.
\end{align}
Note that the latter is a typical situation when $ q_\tau <\infty ;$ this is so, for instance, for gamma-distributed~$\tau$. It turns out that, in this situation, the desired~$\beta$ may or may not exist depending on the distribution of~$\Theta$, .

\begin{figure}[ht]
	\centering
	\includegraphics[scale=.55]{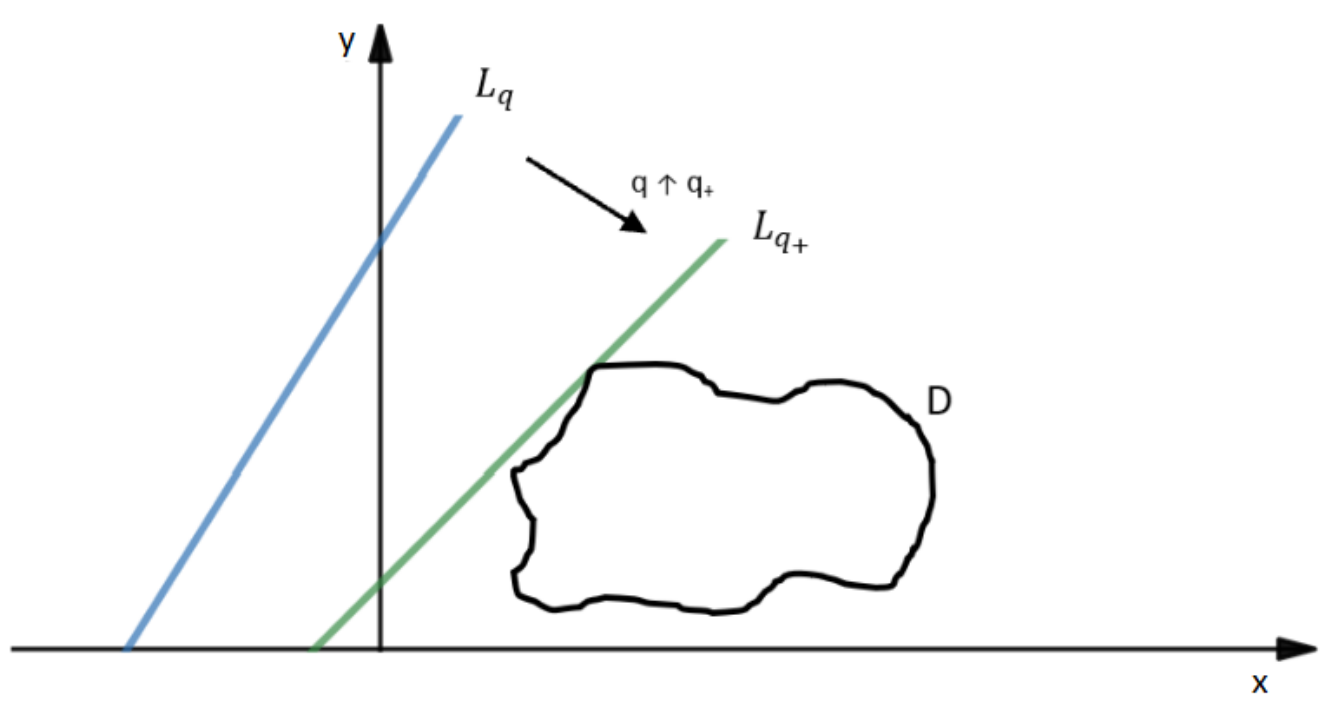}
	%
	%
	\caption{ As $q\uparrow q_+,$ the line $L_q$ approaches the set~$D.$}
	\label{fig:1}       
\end{figure}

Introduce rays $L_q:= \{(x,y)\in \R\times\R^+: \langle u(q), (x,y)\rangle = q_\tau \}$, $q>0.$ Clearly, $(x,y)\in L_q$ iff
\begin{align}
\label{L_eq}
y= \frac{x}{q+1} +\frac{q_\tau}{q(q+1)}, \quad x\ge -\frac{q_\tau }q
\end{align}
(the last inequality is equivalent to $y\ge 0$). Denote by~$D\subseteq [\underline{\mu},\infty)\times [0,\overline{\si}^2/2] $ the support of $\Theta$  and put
\[
q_+:=\inf\{q>0:   L_q \cap D\neq \varnothing\}  .
\]
Note that, as $q$ increases, the ray $L_q$ ``moves" to the right and ``rotates" in the clock-wise direction, and as~$D$ is bounded from the left and from the top, $q_+$ is a finite positive number (see Fig.~\ref{fig:1}: $q_+$ is the value of~$q$ for which~$L_q$ first ``touches"~$D$).

Note that if $q>q_+$ then $\pr (\langle u(q), \Theta \rangle >q_\tau)>0$, so that $\phi_\nu (q)=\infty$ by~\eqref{phi_now}. As for  $q<q_+$ one clearly has $\pr (\langle u(q), \Theta \rangle <q_\tau-\ep)=1$ for some $\ep>0,$ we get $q<q_\nu $ (again by~\eqref{phi_now}). We conclude that $q_+= q_\nu.$

Clearly, $\phi_\nu(q_\nu)=\infty$ is   sufficient for the existence of the Lundberg exponential under the condition that $\exn \nu< 0.$  In view of our assumption~\eqref{tau_infi}, representation~\eqref{phi_now} suggests that whether~$\phi_\nu(q_\nu)$ is infinite or not depends on how strongly the distribution of~$\Theta$ is concentrated in vicinity of the ray~$L_{q_+}.$
To capture this, we introduce the random variable $H$ by setting, for any $\theta_+\in L_{q_+}$,
\[
H:= \langle u(q_+), \theta_+ -\Theta  \rangle
 = q_\tau + q_+\mu -  q_+( q_++1)\si^2/2,
\]
where the second equality was obtained choosing $\theta_+=\theta_0:= (-q_\tau/q_+, 0)\in L_{q_+}$, and denote by~$F_H$ its distribution function. We see that $H\ge 0$ a.s.\ (as the point~$\Theta$ is below the ray~$L_{q_+}$ given by~\eqref{L_eq}) with $q=q_+,$ the value of~$H$ being equal to the Euclidean length of the vector $u(q_+)$ times the distance from~$\Theta$ to~$L_{q_+}$

Now, from~\eqref{phi_now},
\begin{align}
\phi_\nu (q) & =  \exn \phi_\tau (\langle u(q)- u(q_+), \Theta \rangle
 +\langle u(q_+), \Theta -\theta_0 \rangle + \langle u(q_+),\theta_0\rangle)
\notag\\
& =  \exn \phi_\tau (q_\tau -H - \langle u(q_+)- u(q), \Theta \rangle)
=  \exn \phi_\tau (q_\tau -H - \ep \langle a(\ep) , \Theta \rangle),
\label{ff}
\end{align}
where we first noted that $ u(q_+)- u(q)=(q_+-q)(-1, q_+ +q+1) $ and then  put $\ep:=q_+-q,$ $a(\ep):= (-1, 2q_+ +1 -\ep)\to(-1, 2q_+ +1  ) $ as $\ep\downarrow 0.$

There is no monotone dependence on~$\ep$ in the integrand on the right-hand side on~\eqref{ff}, so we need an argument   establishing convergence of these expectations as $\ep\downarrow 0.$ Let $D'':= \{\theta=(x,y)\in D: y<x/(q_++1)\}$. Clearly, $\langle u(q_+), \theta\rangle = -q_+ x + q_+ (q_++1)y <0$ for $\theta\in D'',$ so that $\phi_\tau (\langle u(q_+), \theta\rangle)\le 1$ in that domain  and hence
\[
\exn (\phi_\tau (\langle u(q), \Theta \rangle); \Theta\in D'') \to 
 \exn (\phi_\tau (q_\tau-H ); \Theta\in D''),
\quad   q\uparrow q_+,
\]
by the dominated convergence theorem. Turning to $D':=D\setminus D''$, one can easily verify that there exist $r_{\pm} \in \R$ such that $r_- \le \langle a(\ep) , \theta \rangle \le r_+ $ for all $\theta\in D',$ $\ep\in (0,1)$. Since $\phi_\tau$ is an increasing function, we get
\begin{align*}
 \exn (\phi_\tau (q_\tau -H - r_+\ep  ) ; \Theta\in D')
 & \le  \exn (\phi_\tau (q_\tau -H - \ep \langle a(\ep) , \Theta \rangle); \Theta\in D')
 \\
 &  \le
 \exn (\phi_\tau (q_\tau -H - r_-\ep  ) ; \Theta\in D').
\end{align*}
Now we can apply the monotone convergence theorem to both lower and upper bounds in the last displayed formula  since the integrands in them have monotone dependence on~$\ep.$ We conclude that
\[
\phi_\nu (q)\to \int_0^\infty \phi_\tau (q_\tau - h)dF_H(h), \quad q\uparrow q_+.
\]
Since clearly $\int_{\de}^\infty \phi_\tau (q_\tau - h)dF_H(h)\le \phi_\tau (q_\tau -\de)<\infty$ for any $\de>0,$ we arrive at the following result.

\begin{theorem}\label{T2}
Under the assumptions stated at the beginning of this section,  assume that~\eqref{mu_si_0} and~\eqref{tau_infi} hold true. Then $q_\nu =q_+ $ and, moreover,  $\phi_\nu (q_\nu)=\infty$ iff
 	\[
 	\int_0^\de \phi_\tau (q_\tau - h)dF_H(h)=\infty
 	\]
 	for some $($and then for any$)$ $\de>0$.
\end{theorem}

Thus there must be significant presence of probability mass in vicinity of the tangent to~$D$ line~$L_{q_+}$ to ensure that $\phi_\nu (q_\nu)=\infty.$

It is not hard to get closed-form expressions for~$\phi_\nu$ in several tractable examples in the special case of the  Poisson arrival process with rate~1, which means that $\phi_\tau (q) = 1/(1-q),$ $q<q_\tau:=1 $ (so that~\eqref{tau_infi} is true). In one such example one has   
$\pr (\Theta= (1/j, 1-1/j))=j^{-p}/\zeta (p)$, $j\ge 1,$ for a fixed $p\in \N,$ where~$\zeta$ is the Euler--Riemann zeta function. In this case,  $q_+=(\sqrt{5}-1)/2$ and $D\cap L_{1}=\{(0,1)\}, $ and it turns out that   $\phi_\nu(q_+)=\infty$ iff $p=2,$ in obvious agreement with the claim of Theorem~\ref{T2}. If, further,  one assumes that~$\Theta$ is uniformly distributed in a unit  square~$D$  with vertices at the points $(i,j),$ $i,j\in \{0,1\}$, then again  $q_+=(\sqrt{5}-1)/2$, $D\cap L_{1}$ consists of the single point $(0,1)$ (the vertex of our~$D$ at which  it touches the line~$L_{q_+}$), and one can also derive a closed form expression for~$\phi_\nu$ yielding $\phi_\nu(q_+)<\infty. $ If, however, we rotate the square in the  anticlockwise direction around the  vertex $(0,1)$ until its upper edge runs along the line $L_{q_+}$ (with clearly  the same value of $q_+$ as in the previous examples)  then  one would have $\phi_\nu(q_+)=\infty$, also in agreement with  Theorem~\ref{T2}. In the latter case, there is ``too much probability'' in vicinity of~$L_{q_+}$ (the probability mass in the $\ep$-neighbourhood of that line is $\asymp\ep$ as $\ep\downarrow 0$ compared to $\asymp\ep^2$  in the former case).

%
%

\end{document}